\documentclass[12pt,a4paper]{article}

\usepackage{latexsym}
\usepackage[T1]{fontenc}         
\usepackage[latin1]{inputenc}    
\usepackage[english]{babel}      
\usepackage{indentfirst}         
\usepackage{fancybox}
\usepackage{makeidx}
\usepackage[dvips,twoside,top=20mm,bottom=20mm,left=25mm,right=25mm]{geometry}
\usepackage{bbold}
\usepackage{multirow}
\usepackage{hyperref}            
\usepackage{lmodern}
\usepackage{amsmath, amssymb}    
\usepackage{verbatim}            
\usepackage{fancyhdr}
\usepackage{graphicx}            
\usepackage{subfig}              
\usepackage{titlesec}
\usepackage{stackrel}
\usepackage{cite}
\usepackage[svgnames]{xcolor}

\usepackage{cite}
\usepackage{stackrel}
\usepackage{flafter}

\begin{document}

\title{Forecasting of time data with using fractional Brownian motion\thanks{%
Mathematics Subject Classifications: 60G22, 62M10,91B84.}}
\author{Valeria Bondarenko\thanks{%
Ecole Centrale de Nantes}\ , Victor Bondarenko\thanks{%
National Technical University of Ukraine Kiev Polytechnic University}\ , Kiryl Truskovsky\thanks{%
National Technical University of Ukraine Kiev Polytechnic University}\, Ina Taralova\thanks{%
Ecole Centrale de Nantes }}
\maketitle

\date{\today}
\maketitle

\begin{abstract}
We investigated the quality of forecasting of fractional Brownian motion, and new method for estimating of Hurst exponent is validated.  Stochastic model of the  time series in the form of converted fractional Brownian motion is proposed. The method of checking the adequacy of the proposed model is developed and  short-term forecasting for temporary data is constructed. The research results are implemented in software tools for analysis and modeling of time series.
\end{abstract}


\section{Introduction}

We assume that an observed trajectory $x\left(t\right)$, $0\leq t\leq T$, is an element of  the ensemble of trajectories or an element of function space in the construction of statistical mathematical model. If $x\left(\cdot\right)$ is assumed to be continuous, then this space can be considered a set $C\left(0; T\right)$, which are continuous functions in $\left(0; T\right)$. In other words,
\begin{equation}
x\left(t\right) = \Phi\left(X\left(\cdot\right)\right) \left(t\right)\;,
\label{eq:1}
\end{equation}
where $X\left(s\right)$ is a realization of some random process $\xi\left(s\right)$ with known characteristics, $\Phi$  is a reversible conversion in $C\left(0; T\right)$. $\left(\Phi, \xi\right)$ is called a model of observed data. Process $\xi\left(t\right)$ is basic in model for $x\left(t\right)$. For discrete observation $x_1, \ldots, x_n$ (time series), and in assumption about automodeling  $\xi\left(t\right)$,
$$
x_k = \Phi\left(X\left(\cdot\right)\right) \left(\frac{k}{n}\right)\;,\quad k = 1,\ldots, n\;;\quad\quad \Phi: \mathbb{R}^n \rightarrow \mathbb{R}^n\;.
$$
For the highly oscillating trajectory $x\left(t\right)$, basic process $\xi\left(t\right)$ with unlimited variation is selected. In particular, $\xi\left(t\right) = \sigma B_H\left(t\right)$, where $B_H\left(t\right)$ is a fractional Brownian motion  (fBm),  which was first introduced by B. Mandelbrot in \cite{15, 16} and is defined as a Gaussian random process with zero mean and covariance function:
$$
R\left(t, s\right) = \mathbb{E} B_H\left(t\right)B_H\left(s\right)=\frac{1}{2}\left(t^{2H} + s^{2H} - \left|t - s\right|^{2H}\right)\;,\quad 0 < H < 1\;.
$$

The $n$-dimensional density distribution of fractional Brownian motion looks as follows:
$$
p\left(t_1, \ldots, t_n, x_1, \ldots, x_n\right) = c \text{exp}\left\{-\frac{1}{2}\sum r^{jk}x_j x_k\right\}\;,\quad r_{jk} = R\left(t_j, t_k\right)\;.
$$

Parameter $H\in \left(0; 1\right)$ is called the Hurst exponent of fBm, and the transformation $\Phi^{-1}$ consists of actions that transform fBm realization of the observed trajectory. Using of fractional Brownian motion as a basic process $\xi\left(t\right)$ in the model (\ref{eq:1}) is justified by non-markovian $B_H\left(t\right)$. The first motivation for studies of this process and its applications are considered in \cite{2, 17, 18}. The results of studies of the properties of fractional Brownian motion and its application in models of natural and economic processes are covered in \cite{3,4,12,13,14,26,27,28}. Let's note the reviews \cite{19,25}. 
Research statistics of fractional Brownian motion are quoted below.

Let's choose the model of fractional Brownian  motion for observed time series  $x_1, \ldots, x_n$:
\begin{equation}
x_k = \Phi\left(\left(\sigma B_H\right)\right)\left(\frac{k}{n}\right)\;,
\label{eq:2}
\end{equation}
and transformation $\Phi$ is defined. Let's calculate Hurst exponent of observed time series as $H$ of the basic process $B_H\left(t\right)$. Note that this value depends on transformation $\Phi$. The criteria for adequacy of the representation (\ref{eq:2}) are shown in \cite{3}. From empirical considerations follows that the model (\ref{eq:2}) is suitable for describing the random time data and apriority isn't satisfactory for the approximation of deterministic chaotic sequences. As a rule, the deterministic and stochastic components can be present in observed data.

In present work, a new method of estimation parameters $\sigma$ and $H$ is justified, and the quality of the forecast  is investigated for the observed realization of fractional Brownian motion $x_k = \sigma B_H\left(\frac{k}{n}\right)$. For the real time series, a model is proposed, which uses fractional Brownian motion as a basic process. The criteria of adequacy of this model are developed and short-term forecasting is constructed.

\section{Statistics of Fractional Brownian Motion}
\subsection{Estimation of Parameters}

Let's consider the increments $\xi_k = \sigma\left(B_H\left(\frac{k}{n}\right)- B_H\left(\frac{k-1}{n}\right)\right)$, which form the Gaussian stationary sequence with zero mean and  the correlation matrix   $V = \displaystyle\frac{\sigma^2}{n^{2H}}S$, and elements $s_{jk}$ of the matrix that look as follows:
\begin{equation}
s_{jk} = \rho\left(\xi_j, \xi_k\right) = \frac{1}{2}\left(\left|k - j + 1\right|^{2H} + \left|k -j - 1\right|^{2H} - 2\left|k - j\right|^{2H}\right)\;.
\label{eq:3}
\end{equation}

In particular, the coefficient of correlation between neighbor increments     is
$$
\rho\left(\xi_k, \xi_{k+1}\right) \equiv \rho_1 = 2^{2H - 1} - 1\;.
$$

The limit theorems for sequence $\xi_1, \ldots, \xi_n$ were first proved by Peltier \cite{29}: for statistics
$$
R_{jn} = \frac{1}{n}\sum_{k = 1}^n \left|\xi_k\right|^j\;, \quad j\in \mathbb{N}\;,\quad\quad\quad E_n\left(j\right) = \mathbb{E}R_{jn}=\frac{\sigma^j}{n^{jH}}\frac{2^{\frac{j}{2}} \Gamma\left(\frac{j + 1}{2}\right)}{\sqrt{\pi}}\;,
$$
with probability  1 $\displaystyle\frac{R_{jn}}{E_n\left(j\right)}\rightarrow 1$, $n \rightarrow \infty$.

From the last equation, consistency estimates of parameters  $H$ and $\sigma$ follow:
\begin{equation}
\begin{array}{c}
\hat{H}_n = \displaystyle\frac{\ln \left(\sqrt{\frac{2}{\pi}}\frac{\sigma}{R_{1n}}\right)}{\ln n}\;, \text{with known }\sigma\;,\\
\hat{\sigma}_{1n} = n^H\sqrt{\displaystyle\frac{\pi}{2}}rR_{1n} = 1.25 n^H R_{1n}\;, \text{with known } H\;.
\end{array}
\label{eq:4}
\end{equation}

Let's propose new estimation method of fractional Brownian motion, by observed data $x_1, \ldots, x_n$, two unknown parameters  $\sigma$, $H$.

Let's assume:
$$
\begin{array}{c}
y_1, \ldots, y_n \text{ are the increments}\;; \quad y_k = x_k - x_{k - 1}\;,\\
Q\left(H\right) = \displaystyle\frac{0.8}{R_{1n}}\sqrt{\displaystyle\frac{\left(S^{-1}y, y\right)}{n}}\;,
\end{array}
$$
where matrix $S \equiv S_H$ is defined by (\ref{eq:3}), $y$ is a vector of increments.

\begin{center}
\textbf{Statement}
\end{center}

Statistic
\begin{equation}
\hat{H} = \arg\min\left|Q\left(H\right) - 1\right|
\label{eq:5}
\end{equation}
is a consistent estimator of the parameter $H$.

\begin{center}
\textbf{Proof}
\end{center}

$\varepsilon$ is the canonical Gaussian vector with the following characteristics:
$$
\mathbb{E} \varepsilon = 0\;,\quad \mathbb{E}\left(\varepsilon,u\right) \left(\varepsilon, \upsilon\right) = \left(u, \upsilon\right)\;,\quad \dim \varepsilon = n\;.
$$
Then, $\mathbf{y} = V^{\frac{1}{2}} \mathbf{\varepsilon}$, therefore
$$
n = \mathbb{E}\left(\mathbf{\varepsilon},\mathbf{\varepsilon}\right) = \mathbb{E}\left(V^{-1}\mathbf{y},\mathbf{y}\right) = \displaystyle\frac{n^{2H}}{\sigma^2}\mathbb{E}\left(S^{-1}\mathbf{y},\mathbf{y}\right)\;.
$$
And consequently the statistic
$$
\hat{\sigma}_{2n}^{2} = \left(n\right)^{2H - 1}\left(S^{-1}\mathbf{y},\mathbf{y}\right)\;,
$$
and here statistic $\left(n\right)^{2H - 1}\left(S^{-1}y, y\right)$ is an unbiased estimate of the parameter $\sigma^2$. The dispersion of estimate
\begin{equation}
\hat{\sigma_{2n}} = \sqrt{n^{2H - 1}\left(S^{-1}\mathbf{y},\mathbf{y}\right)}\;.
\label{eq:6}
\end{equation}
is calculated using the formula of integration by parts (\cite[p. 206]{5}). From (\ref{eq:4}) and (\ref{eq:6}), it follows that
$$
\frac{\hat{\sigma}_{2n}}{\hat{\sigma}_{1n}} = \frac{0.8}{R_{1n}}\sqrt{\frac{\left(S^{-1}\mathbf{y}, \mathbf{y}\right)}{n}} = Q\left(H\right)\;,
$$
and consistency of estimates means that $\lim_n Q\left(H\right) = 1$, where $H$ is a Hurst exponent of observed fractional Brownian motion. $\blacksquare$

The implementation of the corresponding algorithm is to choose such a value of argument $H$ in $Q\left(H\right)$, where $\left|Q\left(H\right) - 1\right| \rightarrow \min$.

The efficiency of the algorithm is confirmed by numerical experiment. The statistical values are shown in Table \ref{table:1}.
$$
q_{kj} = \frac{0.8}{R_{1n}} \sqrt{\frac{\left(S^{-1}_j\mathbf{z}_k, \mathbf{z}_k\right)}{n}}\;,
$$
where $\mathbf{z}_k$ is a generated vector of increments fBm with Hurst exponent $H_k$, $S_j$ is the normalized correlation matrix, corresponding to the index fBm with Hurst exponent $H_j$. For every $H_k$, values $q_{kj}$ are calculated with the selection of parameter $H_j$ with step $\Delta H_j = 0.1$. Generation $\mathbf{z}_k$ is performed with the following parameters:
$$
n = 200\;; n=1000\;; H_k = 0.1;0.3;0.7;0.9\;.
$$

\begin{table}[!htp]
\centering
\caption{Efficiency of evaluation method}
\begin{tabular}{|c|c|c|c|c|c|c|c|c|c|}
\hline
\multicolumn{2}{|r|}{$H_j$} & $0.1$ & $0.2$ & $0.3$ & $0.4$ & $0.6$ & $0.7$ & $0.8$ & $0.9$ \\
\multicolumn{2}{|l|}{$H_k$} & & & & & & & & \\
\hline
$0.1$ & $n = 200$ & $1.004$ & $0.95$ & $0.93$ & $0.91$ & $1.11$ & $1.21$ & $1.45$ & $2.03$\\
 & $n = 1000$ & $1.002$ & $0.95$ & $0.93$ & $0.91$ & $1.10$ & $1.22$ & $1.47$ & $2.06$\\
\hline
$0.3$ & $n = 200$ & $1.29$ & $1.07$ & $0.98$ & $0.95$ & $1.07$ & $1.16$ & $1.38$ & $1.92$\\
 & $n = 1000$ & $1.28$ & $1.07$ & $0.99$ & $0.94$ & $1.08$ & $1.17$ & $1.40$ & $1.94$\\
\hline
$0.7$ & $n = 200$ & $4.06$ & $2.25$ & $1.52$ & $1.18$ & $0.94$ & $0.97$ & $1.09$ & $1.44$\\
 & $n = 1000$ & $7.04$ & $3.19$ & $1.79$ & $1.22$ & $0.92$ & $0.98$ & $1.08$ & $1.43$\\
\hline
$0.9$ & $n = 200$ & $7.67$ & $3.88$ & $2.26$ & $1.43$ & $0.72$ & $0.74$ & $0.75$ & $0.97$\\
 & $n = 1000$ & $9.10$ & $4.13$ & $2.24$ & $1.40$ & $0.77$ & $0.78$ & $0.83$ & $1.07$\\
\hline
\end{tabular}
\label{table:1}
\end{table}

Analysis of  data in Table \ref{table:1} shows that, for each $H_k$ (in the fixed line),
$$
\left|q_k - 1\right| \rightarrow \min\;, \text{if } H_j = H_k\text{, so }\hat{H_k} = H_j\;.
$$

Note: In the works of J.-F. Coeurjolly \cite{9,10,11}, supplemented by the work \cite{1}, another method of estimating Hurst exponent is justified, embedded in Package \emph{dvfBm} (https://cran.r-project.org/web/packages/dvfBm/dvfBm.pdf). Let's denote by $\hat{H}_1$ an estimate of the proposed method in this paper and  by $\hat{H}_2$ an estimate by J.-F. Coeurjolly. The comparison of these estimates shows that their deviation is not more than 5\%.

\subsection{Forecast of Fractional Brownian Motion}

Let's observe the trajectory of a random process $x\left(t\right)$, $0 \leq t \leq T$. The random value $\hat{X}\left(T + \tau\right)$ is called an optimal forecast of process in point  $T + \tau$, if
$$
\mathbb{E}\left(\hat{X}\left(T + \tau\right) - X\left(T + \tau\right)\right)^2 = \stackrel[\xi]{}{\min}\mathbb{E}\left(\xi - X\left(T + \tau\right)\right)^2\;.
$$

The optimal forecast is defined by the formula of conditional mean:
\begin{equation}
\hat{X}\left(T + \tau\right) = \mathbb{E}\left(X\left(T + \tau\right)|X\left(t\right)\;,\quad 0 \leq t \leq T\right)\;.
\label{eq:7}
\end{equation}

In some cases, (\ref{eq:7}) assumes an explicit expression. Let's consider Gaussian random vector $\xi=\left(\xi_1, \ldots, \xi_n\right)$, $\xi=\left(\eta, \sigma\right)$;, $\text{dim}\eta = m$, $\text{dim}\sigma = n - m$, $\xi\sim \aleph\left(0; S\right)$, $\eta\sim\aleph\left(0;A\right)$, $\xi\sim\aleph\left(0;D\right)$, so correlation operator of vector  $\xi$ is a block matrix
$$
S = \left[\begin{array}{cc}
A & B\\
C & D
\end{array}\right]\;,
$$
where matrix elements  represent the cross-correlation of coordinates $\eta$ and $\sigma$. If $\eta$ is observed  and  $\sigma$ is estimated vector, then optimal forecast coincides with the linear estimation, and the equation (\ref{eq:7}) takes the following form:
$$
\hat{\sigma} = \mathbb{E}\left(\sigma|\eta\right)=CA^{-1}\eta\;,
$$
or in the coordinate form:
\begin{equation}
\hat{\xi}_{m+j} = \sum_{k = 1}^m\sum_{i = 1}^m s_{m+j,k}a^{ki}\xi_i\;,\quad j = 1,\ldots,n-m\;.
\label{eq:8}
\end{equation}

For one-step extrapolation, $m = n - 1$, $B$ is a vector, convector $C = B^T$, $D = \mathbb{E}\xi_n^2$, and mean-square absolute error of forecast $\delta$ is defined by the following formula:
$$
\delta^2 = \mathbb{E}\left(\hat{\xi}_n - \xi_n\right)^2 = D - \left(A^{-1}B, B\right)\;,
$$
and $\delta D^{\beta 0.5}$ is an error.

We can construct the forecast of fBm  for its increments
$$
\xi_k = y_k = B_H\left(\frac{k}{n}\right) - B_H\left(\frac{k - 1}{n}\right)\;,
$$
as well as for the values of fractional Brownian motion: $\xi_k = B_H\left(\frac{k}{n}\right)$.

In the first case, the elements of the matrix $S$ are defined by (\ref{eq:3}), and formula (\ref{eq:8}) 
takes the following form:
\begin{equation}
\begin{array}{c}
\hat{y}_{m+j} = \displaystyle\sum_{k = 1}^m\displaystyle\sum_{i = 1}^m \left(\displaystyle\frac{\left(m+j-k+1\right)^{2H} + \left(m+j-k-1\right)^{2H}}{2} - \left(m + j - k\right)^{2H}\right)a^{ki}y_i\;,\\
\quad j = 1,\ldots,r\;,\quad r = m-n\;.
\end{array}
\label{eq:9}
\end{equation}

The elements $s_{jk}$ of correlation matrix  $S$ in extrapolating the values of fractional Brownian motion are defined by the following equation:
\begin{equation}
s_{jk} = 0.5\left(j^{2H} + k^{2H} - \left(k - j\right)^{2H}\right)\;.
\label{eq:10}
\end{equation}

The numerical experiment was performed with the simulated data for determining the quality of forecast. Forecast was constructed for 8 steps by the learning sample.

The results of forecast $\left\{y_k\right\}$ by the formula (\ref{eq:9}) are not satisfactory: the absolute error $\delta_j = \left| \frac{\hat{y}_{m+j} - y_{m+j}}{y_{m+j}}\right|$, $j = 1, \ldots, 8$, equals $0.8$--$1.2$, and doesn't depend on the size of the learning sample.

The calculation of forecasting values $\hat{x}_{m+j}$ by formula (\ref{eq:8}) with the matrix defined by equation (\ref{eq:10}) leads to the following expected result. The forecast of antipersistent process ($H < 0.5$) is not satisfactory, the error of forecast does not depend on the size of the learning sample. The quality of forecast improves with increasing $m$ for the persistent process. Appropriate data are given in Table \ref{table:2}, which shows  the values of the  relative error $\delta_j = \left| \frac{\hat{x}_{m+j} - x_{m+j}}{x_{m+j}}\right|$, $j = 1, \ldots, 8$, for $H = 0.3$, $H = 0.7$, $H = 0.9$, $m = 100$, $m = 500$, $m = 1000$.

\begin{table}[!htp]
\centering
\caption{}
\begin{tabular}{|c|c|c|c|c|c|c|c|c|c|}
\hline
 \multicolumn{2}{|l|}{$H$} & 1 & 2 & 3 & 4 & 5 & 6 & 7 & 8 \\
\hline
$0.3$ & $m = 100$ & $0.67$ & $0.71$ & $0.45$ & $0.07$ & $2.32$ & $1.49$ & $1.58$ & $0.11$\\
 & $m = 500$ & $0.25$ & $0.16$ & $0.06$ & $0.17$ & $1.68$ & $0.61$ & $0.71$ & $0.69$\\
 & $m = 1000$ & $0.04$ & $0.12$ & $0.16$ & $0.01$ & $0.10$ & $0.04$ & $0.18$ & $0.17$\\
\hline
$0.7$ & $m = 100$ & $0.02$ & $0.08$ & $0.17$ & $0.27$ & $0.29$ & $0.35$ & $0.50$ & $0.58$\\
 & $m = 500$ & $0.015$ & $0.012$ & $0.009$ & $0.003$ & $0.027$ & $0.001$ & $0.005$ & $0.007$\\
 & $m = 1000$ & $0.008$ & $0.015$ & $0.031$ & $0.017$ & $0.006$ & $0.001$ & $0.021$ & $0.017$\\
\hline
$0.9$ & $m = 100$ & $0.02$ & $0.07$ & $0.10$ & $0.18$ & $0.21$ & $0.24$ & $0.26$ & $0.31$\\
 & $m = 500$ & $0.001$ & $0.001$ & $0.001$ & $0.001$ & $0.01$ & $0.01$ & $0.01$ & $0.01$\\
 & $m = 1000$ & $0.001$ & $0.01$ & $0.02$ & $0.04$ & $0.05$ & $0.05$ & $0.07$ & $0.07$\\
\hline
\end{tabular}
\label{table:2}
\end{table}

\textbf{Note:} The data of forecast for  $H = 0.9$, $m = 1000$ may contain some errors due to poor conditioning of matrix $S_H$: determinant of this matrix is a decreasing function  of $m$ and $H$ for $H > 0.5$. So, for $m = 500$, $\text{det} S_{0.9} \sim 10^{-190}$.

\section{Limit Theorems and Applications}
\label{sec:limit}

Let $B\left(t\right)$, $0 \leq t \leq 1$ be a fractional Brownian motion with Hurst  exponent $H$. Let's consider the normalized increments
$$
\xi_k = n^H\left(B\left(\frac{k}{n}\right) - B\left(\frac{k - 1} {n}\right)\right)\sim\aleph\left(0;1\right)\;.
$$

In the series of papers \cite{7,8,20,21,22,23,24},  some limit theorems for the functions of these increments were proven. Let's denote
$$
\alpha_k = n^H B\left(\frac{k}{n}\right) = \sum_{j = 1}^{k - 1}\xi_j\;.
$$
There is a mean-square convergence:
\begin{equation}
\begin{array}{c}
\displaystyle\frac{1}{n}\displaystyle\sum_{k = 1}^n \alpha_k \xi_k^3 \to -\displaystyle\frac{3}{2}\;, \quad H \in \left(0; \displaystyle\frac{1}{2}\right)\;,\\
\displaystyle\frac{1}{n^{1 + H}}\displaystyle\sum\alpha_k^2\xi_k^3 \to 3\eta\;, \quad H \in \left(0;\displaystyle\frac{1}{2}\right)\;,
\end{array}
\label{eq:11}
\end{equation}
where
$$
\eta \sim \aleph\left(0;\frac{1}{2H + 2}\right)\;;
$$
$$
\frac{1}{n^{2H}} \sum_{k = 1}^n \alpha_k \xi_k^3 \to \frac{3}{2}B^2\left(1\right)\;, \quad H \in \left(\frac{1}{2};1\right)\;.
$$

These limit relations allow us to check the statistical hypothesis $T = \left\{\text{the investigated}\right.$ $\left.\text{time series }x_1, \ldots, x_n\text{ is an implementation of fBm}\right\}$.

\textbf{The algorithm of checking is as follows} (with known $H$)

Let's consider the increments $y_k = x_k - x_{k - 1}$, statistics $R_{1n}\left(y\right) = \displaystyle\frac{1}{n} \displaystyle\sum_{k = 1}^n\left|y_k\right|$, and estimate $\hat{\sigma}$ by formula (\ref{eq:4}).

Let's normalize the increments and assume:
$$
z_k = \left(\hat{\sigma}\right)^{-1}n^H y_k = \frac{0.8}{R_{1n}}y_k\;.
$$

We assume that the hypothesis $T$ holds:
\begin{equation}
z_k = \xi_k = n^H\left(B\left(\frac{k}{n}\right) - B\left(\frac{k - 1}{n}\right)\right)\;.
\label{eq:12}
\end{equation}

Assume $v_k = \displaystyle\sum_{j = 1}^{k - 1}z_j$ and calculate the statistics
\begin{equation}
\begin{array}{cc}
A_n = \displaystyle\frac{1}{n}\displaystyle\sum v_kz_k^3\;, & H\in \left(0;\frac{1}{2}\right)\;;\\
B_n = \displaystyle\frac{1}{n^{1 + H}}\displaystyle\sum v_k^2z_k^3\;, & H\in \left(0;\frac{1}{2}\right)\;;\\
D_n = \displaystyle\frac{1}{n^{2H}}\displaystyle\sum v_k z_k^3\;, & H\in \left(\frac{1}{2};1\right)\;.
\end{array}
\label{eq:13}
\end{equation}

If hypothesis $T$ is true, then there is convergence:
$$
A_n \to -1.5\;; \quad 
B_n \to 3\eta\;; \quad
D_n \to \frac{3}{2}B^2\left(1\right)\;.
$$
The decision about the hypothesis $T$ is taken  by comparing the actual  values of  statistics with their limiting theoretical values. Let's determine the deviation from the limit value $\delta = \left|A_n + 1.5\right|$ for statistic $A_n$; the limit distribution functions for statistics $B_n$, $D_n$:
$$
F_1\left(x\right) = P\left\{3\eta < x\right\} = \Phi\left(\frac{x}{3d}\right)\;, \quad F_2\left(x\right) = 2\Phi\left(\sqrt{\frac{2}{3}}x\right) - 1\;, x > 0\;,
$$
where $\Phi$ is Laplace function, $d = \left(2H + 2\right)^{-0.5}$.

Hypothesis $T$ is accepted, if
\begin{equation}
\delta < \beta_0\;, \left|B_n\right| < \beta_1\;, H < 0.5\;; \quad 0 < D_n < \beta_2\;, H > 0.5\;,
\label{eq:14}
\end{equation}
where $\beta_1$, $\beta_2$ are quantiles of distributions of $F_1$, $F_2$,  corresponding to the selected level of significance  $\alpha = 0.1$. Then,
$$
\beta_1 = \frac{4.95}{\sqrt{2H + 2}}\;, \quad \beta_2 = 4.08\;.
$$

The rate of convergence of statistics to the limit has been tested by numerical experiment  for the first example (``ideal case''):
$$
z_k = \left(\hat{\sigma}^{-1}\right)n^H\left(X\left(k\right) - X\left(k - 1\right)\right)\;; \quad X\left(t\right) = \sigma B_H\left(t\right)\;,
$$
where the values of fractional Brownian motion were obtained by simulation. The values of the control statistics $A_n$, $B_n$, $D_n$ are shown in Table \ref{table:3}.

\begin{table}[!htp]
\centering
\caption{Values of control statistics}
\begin{tabular}{|c|c|c|c|c|c|}
\hline
 \multicolumn{2}{|l|}{$H$}  & $A_n$ & $B_n$ & $D_n$ & $\beta_1$\\
\hline
$0.1$ & $n = 200$ & $-1.30$ & $0.84$ & & $3.34$\\
 & $n = 1000$ & $-1.32$ & $2.63$ &  & $3.34$\\
\hline
$0.2$ & $n = 200$ & $-1.21$ & $0.81$ & & $3.20$\\
 & $n = 1000$ & $-1.35$ & $1.74$ &  & $3.20$\\
\hline
$0.3$ & $n = 200$ & $-2.00$ & $0.37$ &  & $3.07$\\
 & $n = 1000$ &  $-1.10$ & $0.50$ &  & $3.07$\\
\hline
$0.4$ & $n = 200$ & $-0.55$ & $1.26$ &  & $2.96$\\
 & $n = 1000$ & $-2.51$ & $0.83$ &  & $2.96$\\
\hline
$0.6$ & $n = 200$ & & & $1.75$ &\\
 & $n = 1000$ & & & $1.03$ &\\
\hline
$0.7$ & $n = 200$ & & & $1.23$ &\\
 & $n = 1000$ & & & $0.67$ &\\
\hline
$0.8$ & $n = 200$ & & & $1.05$ &\\
 & $n = 1000$ & & & $0.52$ &\\
\hline
$0.9$ & $n = 200$ & & & $0.48$ &\\
 & $n = 1000$ & & & $0.04$ &\\
\hline
\end{tabular}
\label{table:3}
\end{table}

From Table \ref{table:3}, it follows that
$$
\left|B_n\right| < \frac{4.95}{\sqrt{2H + 2}} = \beta_B\;, \quad H < 0.5\;; \quad 0 < D_n < 4.08 = \beta_D\;, \quad H > 0.5\;,
$$
and deviation  $\delta$, $H < 0.5$, is an increasing function of $H$ (for $H = 0.4$, $\delta \approx 1$).

\textbf{The second example} is a deterministic logistic chaotic sequence $x_{k+1} = 4x_k\left(1 - x_k\right)$, $k = 1, \ldots, 1049$. By procedure of estimation, we obtained $\hat{H} = 0.15$. The control statistics are as follows: $A_n = 0.6 > 0$, $\left|B_n\right| = 1.9 > \beta_B = 0.08$. Hypothesis $T$ is rejected.

\textbf{The third example.} Assume that the observed values are an  additive mixture of the deterministic chaotic and random sequences:
$$
x_k = u_k + av_k\;,
$$
where $u_k$ are the values of a dynamical system, $v_k$ are the values of a random process.

Sequences $\left\{u_k\right\}$, $\left\{v_k\right\}$ are normalized by energy, therefore $\frac{1}{n}\sum u_k^2 = \frac{1}{n} \sum v_k^2 = 1$. Then, the value $a$ determines stochastic share in the observed data.

In the example, $u_k = 4u_{k - 1}\left(1 - u_{k - 1}\right)$, $v_k = \sigma B_H\left(\frac{k}{n}\right)$.

The stochastic sequence $v_k$ is generated with  $H_{fBm} = 0.1$--$0.9$. Table \ref{table:4} shows estimate  $\hat{H}$ of mixture and values of control statistics.

\begin{table}[!htp]
\centering
\caption{Control statistics of mixture ($a = 1$, $a = 2$, $n = 2000$)}
\begin{tabular}{|c|c|c|c|c|c|c|}
\hline
 \multicolumn{2}{|l|}{$H$}  & $H$ & $A_n$ & $B_n$ & $D_n$ & $\beta_1$\\
\hline
$0.1$ & $a = 1$ & $0.6$ & $-1.94$ & $-0.07$ & $-0.43$ & $2.77$\\
 & $a = 2$ & $0.1$ & $-1.60$ & $-0.40$ & $-697$ & $3.34$\\
\hline
$0.2$ & $a = 1$ & $0.15$ & $-5.35$ & $-15.3$ & $-1095$ & $3.26$\\
 & $a = 2$ & $0.15$ & $-3.19$ & $-6.17$ & $-652$ & $3.26$\\
\hline
$0.3$ & $a = 1$ & $0.6$ & $-2.54$ & $-0.12$ & $-0.56$ & $2.77$\\
 & $a = 2$ &  $0.2$ & $-5.0$ & $-9.20$ & $-477$ & $3.19$\\
\hline
$0.4$ & $a = 1$ & $0.15$ & $-4.50$ & $-7.38$ & $-920$ & $3.26$\\
 & $a = 2$ & $0.15$ & $-2.33$ & $-1.03$ & $-475$ & $3.26$\\
\hline
$0.6$ & $a = 1$ & $0.6$ & $-0.91$ & $-0.01$ & $-0.20$ & $2.77$\\
 & $a = 2$ & $0.15$ & $-4.0$ & $15.8$ & $-813$ & $3.26$\\
\hline
$0.7$ & $a = 1$ & $0.6$ & $-1.35$ & $-0.03$ & $-0.30$ & $2.77$\\
 & $a = 2$ & $0.1$ & $-1.08$ & $0.47$ & $-470$ & $3.34$\\
\hline
$0.8$ & $a = 1$ & $0.6$ & $-1.37$ & $-0.03$ & $-0.30$ & $2.77$\\
 & $a = 2$ & $0.6$ & $-0.68$ & $-0.01$ & $-0.15$ & $2.77$\\
\hline
$0.9$ & $a = 1$ & $0.6$ & $-1.45$ & $-0.03$ & $-0.32$ & $2.77$\\
 & $a = 2$ & $0.6$ & $-1.91$ & $-0.07$ & $-0.42$ & $2.77$\\
\hline
\end{tabular}
\label{table:4}
\end{table}

The table data show the ``aggressiveness'' of the chaotic component in relation to  stochastic for $H_{fBm} \geq 0.2$. Inequalities (\ref{eq:14}) are not satisfied for these values of fBm and character of the mixture determines the logistic sequence. The deviation of statistics from the limit values is the same as for the ``pure'' fractional Brownian motion (Table \ref{table:3}) (for $H_{fBm} = 0.1$).

Conclusion: Persistence ($\hat{H} > 0.5$) of investigated time series  ($D_n < \beta_2$) means it has stochastic nature; antipersistent ($\hat{H} = 0.1$--$0.2$, $A_n \approx A$, $\left|B_n\right| < \beta_1$) admits the existence of the chaotic component.

\section{The Real Data: Approximation and Forecast}

Construction of the model (\ref{eq:2}) for real-time series $x_0 = 0$, $x_1,\ldots,x_n$, $\overline{x} = 0$ is in choosing the transformation $\Phi$ and checking the adequacy of the model by criterion (\ref{eq:14}). The transformation $\Phi^{-1}$ is defined on the vector of increments
$$
y = \left(y_1, \ldots, y_n\right)\;, \quad y_k = x_k - x_{k - 1}\;, \quad \tilde{y} = \left(\tilde{y}_1, \ldots, \tilde{y}_n\right) = \Phi^{-1}\left(y\right)\;,
$$
where $\tilde{y}_k$ are the increments of fractional Brownian motion. The procedure of constructing the model is called ``algorithm for approximating the time series $s_0 = 0$, $s_1, \ldots, s_n$ by fractional Brownian motion,'' which consists of the following:
\begin{enumerate}
	\item Primary conversion  $\psi$ on initial data $s_0 = 0$, $s_1, \ldots, s_n$, which is leading  new sequence $x_0 = 0$, $x_1,\ldots,x_n$, $\overline{x} = 0$ ($x_k = \psi\left(s_k\right)$), and calculation of the increments $y_k = x_k - x_{k - 1}$. In particular, the transformation  $\psi$ may contain a logarithm and removing approximation of the trend ($S_k > 0$, $x_k = \log S_k - M_k$).
	\item Selection of operator $\Phi^{-1}$, which is converting the increments $\tilde{y}_k$ in new sequence $\left(\tilde{y}_1, \ldots, \tilde{y}_n\right)$:
\begin{equation}
\tilde{y}_k = \sigma\left(B\left(\frac{k}{n} - B\left(\frac{k - 1}{n}\right)\right)\right)\;,
\label{eq:15}
\end{equation}
and construction of the new time series $u_k = \displaystyle\sum_{j = 1}^k \tilde{y}_j$.
\item Estimation of $H$ exponent by (\ref{eq:5}), where $y = \left\{\tilde{y}_k\right\}$.
\item Investigation of the adequacy of the proposed model or checking the statistical hypothesis (\ref{eq:15}). Adequacy is checked by methods described in Section \ref{sec:limit}, which are reduced to the calculation of control statistics (\ref{eq:13}):
$$
z_k = \frac{0.8}{R_{1n}\left(\tilde{y}\right)}\tilde{y}_k\;,
$$
and comparison of these statistics with the limit values. The hypothesis (\ref{eq:15}) is accepted if relations (\ref{eq:14}) hold.
\item Forecast for $r$ steps for converted time series $u_1, \ldots, u_n$, based on this model:
\begin{equation}
\hat{u}_{m+j} = \sum_{k = 1}^m \sum_{i = 1}^m s_{m+j,k}s^{ki}u_i\;, \quad j = 1, \ldots, r\;,
		\tag{8A}\
\label{eq:8a}
\end{equation}
where $m$ is the size of the learning sample, and the elements $s_{jk}$ of correlation matrix $S$ are defined by the equality (\ref{eq:10}). The reverse transition to the forecast $\hat{s}_{m + 1}, \ldots, \hat{s}_{m + r}$ of initial data is performed by the following procedure:
	\item[5a.] Calculation of the vector of increments:
	$$
	\begin{array}{c}
	v_1 = \tilde{u}_{m+1} - u_m\;,\\
	v_2 = \tilde{u}_{m+2} - \tilde{u}_{m+1}\;,\\
	\ldots\\
	v_r = \tilde{u}_{m+r} - \tilde{u}_{m+r-1}\;,
	\end{array}
	$$
and its conversion into a new vector:
$$
w = \left(w_1, \ldots, w_r\right)\;, \quad w = \Phi\left(v\right)\;.
$$
\item[5b.] Construction of forecast of an auxiliary time series:
$$
\tilde{x}_{m+j} = x_m + \sum_{k = 1}^j w_k\;, \quad j = 1, \ldots, r\;,
$$
and initial time series
$$
\left(\hat{s}_{m+1}, \ldots, \hat{s}_{m+r}\right) = \psi^{-1}\left(\hat{x}_{m+1}, \ldots, \hat{x}_{m+r}\right)\;.
$$
\end{enumerate}

It is necessary to investigate the sequence $\left\{y_1, \ldots, y_n\right\}$ for the realization Selection 2.

In \cite{6}, the following method  was proposed for constructing a one-dimensional transformation $\Phi^{-1}$ for the sample $\left\{y_1, \ldots, y_n\right\}$, $n \sim 200$--$1000$. Let's consider the kurtosis:
$$
d\left(y\right) = \frac{R^2_{1n}\left(y\right)}{R_{2n}\left(y\right)}\;.
$$
If $d_n$ is significantly different from $\frac{2}{\pi}$, let's replace the time series $\left\{y_1, \ldots, y_n\right\}$ with the new sequence $\left\{\tilde{y}_1, \ldots, \tilde{y}_n\right\}$ by the following formula:
\begin{equation}
\tilde{y}_k = \text{sgn} y_k \left|y_k\right|^{\frac{1}{\lambda}}\;, \quad y_k = \text{sgn} y_k \left|\tilde{y}_k\right|^{\lambda}\;, \quad \lambda > 0\;,
\label{eq:16}
\end{equation}
where parameter  $\lambda$ is defined from the following equation:
$$
d = \frac{1}{\sqrt{\pi}}\frac{\Gamma^2\left(\frac{\lambda + 1}{2}\right)}{\Gamma\left(\lambda + \frac{1}{2}\right)}\;; \quad d\left(y\right) = \frac{R^2_{1n}\left(y\right)}{R_{2n}\left(y\right)}\approx\frac{2}{\pi}\;.
$$

Thus, the proposed approximation leads to the following  model of original time series:
$$
x_k = \sum_{j = 1}^k \text{sgn}y_j\left|\tilde{y}_j\right|^{\lambda}\;.
$$

\begin{figure}[!htp]%
\centering
\caption{Concentration of carbon dioxide}%
\includegraphics[scale=1.5]{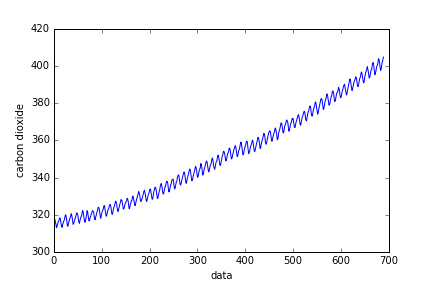}%
\label{fig:1}%
\end{figure}

Let's consider two examples of real data:
\begin{enumerate}
	\item Carbon dioxide (http://climate.nasa.gov/vital-signs/carbon-dioxide/) from 1.\-03.\-1958 to 1.06.2016, 693 data points (Fig. \ref{fig:1}).
	\item The prices of soybean oil (LAMETA, Department of Economics, University of Montpellier)  from 01.01.1960 to 01.09.2011, 610 data points. (Fig. \ref{fig:1})
\end{enumerate}

\begin{figure}[!htp]%
\centering
\caption{The prices of soybean oil, (\$/mt)}%
\includegraphics[scale=0.8]{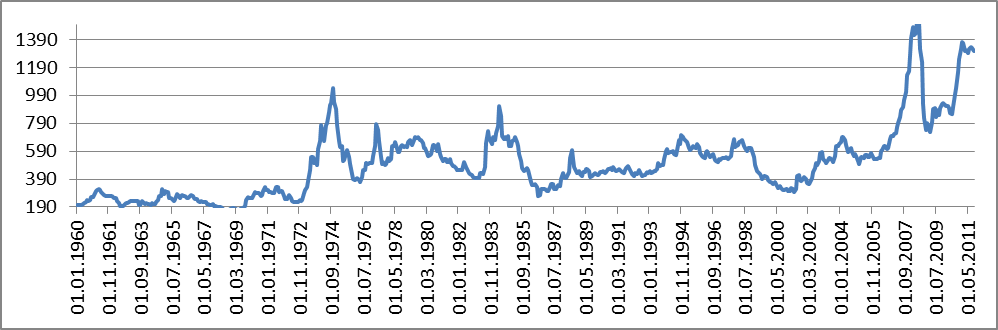}%
\label{fig:2}%
\end{figure}

The initial algorithm of transformation $\psi$ is as follows: the values $\log s_k$ are divided into time windows (a linear approximation of the trend is constructed in every window).

The initial research on stationary increments consists in calculating the correlation coefficient $\hat{p}_1$ for three time windows by the following formula:
$$
\hat{p}_1 = \frac{\sum y_j y_{j + 1}}{\sum y_j^2}\;.
$$

For the first example: $\hat{p}_1 = 0.24$--$0.25$. For the second example: $\hat{p}_1 = 0.63$--$0.65$.

The increments form a stationary sequence, since the values of   $\hat{p}_1$ don't  depend on the number of the window.

The results of calculation are shown in the Table \ref{table:5}.
 
\begin{table}[!htp]
\centering
\caption{Characteristics of real data}
\begin{tabular}{|c|c|c|c|c|c|c|c|c|}
\hline
 Example & Log Mean trend & $d\left(y\right)$ & $\lambda$ & $\hat{H}$ & $A_n$ & $B_n$ & $D_n$ & $\beta_1$\\
\hline
carbon & $0.0037$ & $0.74$ & $0.75$ & $0.75$ & $26$ & $-5.3$ & $1.0$ & $2.66$\\
\hline
Huile & $0.022$ & $0.50$ & $1.37$ & $0.65$ & $0.04$ & $0.14$ & $0.01$ & $2.72$\\
\hline
\end{tabular}
\label{table:5}
\end{table}

The forecast of real data was  performed for the size of learning sample $r = 4$, $m = 200; 400; 600$. The transition to the forecast of initial data is determined by the following formula:
$$
\hat{s}_{m+j} = \text{exp}\left\{\hat{x}_{m+j}+M_{k+j}\right\}\;.
$$

Values of forecast error are shown in Table \ref{table:6}.

\begin{table}[!htp]
\centering
\caption{Values of forecast error $\delta_{m+k}$}
\begin{tabular}{|c|c|c|c|c|c|}
\hline
 \multicolumn{2}{|l|}{$m$} & $k = 1$ & $k = 2$ & $k = 3$ & $k = 4$\\
\hline
$200$ & carbon dioxide & $0.003$ & $0.004$ & $0.002$ & $0.003$\\
\cline{2-6}
 & Huile & $0.09$ & $0.07$ & $0.09$ & $0.09$\\
\hline
$400$ & carbon dioxide & $0.0006$ & $0.002$ & $0.005$ & $0.005$\\
\cline{2-6}
& Huile & $0.009$ & $0.012$ & $0.09$ & $0.04$\\
\hline
$600$ & carbon dioxide & $0.001$ & $0.006$ & $0.01$ & $0.01$\\
\cline{2-6}
& Huile & $0.02$ & $0.04$ & $0.04$ & $0.06$\\
\hline
\end{tabular}
\label{table:6}
\end{table}

Table \ref{table:6} confirms the satisfactory quality of forecast.

\section*{Conclusions}

The proposed model of real time series with fractional Brownian motion as a basic process is effective, if the increments of the observed data have the property of stationarity. Considered examples of physical and financial nature allow  an approximation by a persistent process, which is confirmed by checking the adequacy of  model. Constructed short-term forecast is satisfactory.


\begin{thebibliography}{99}
\bibitem{1}
Achard S. and  Coeurjolly J.-F. Discrete variations of the fractional Brownian in the presence of outliers and an additive noise.  Statistics Surveys (IMS),  Vol. 4, 2009, P. 117--147.

\bibitem{2}
Beran J. Statistics for Long-Memory Processes / Beran J. --- Chapman and Hall. --- 1995. --- 315 p.  

\bibitem{3}
Bezborodov V., Mishura Y., Luca Di Persio. Option pricing with fractional stochastic volatility and discontinuous payoff function of polynomial growt. arXiv: 16.07.07392[math.PR] (Submitted on 25 Jul 2016)

\bibitem{4}
F. Biagini, Y. Hu, B. \O ksendal, T. Zhang. Stochastic calculus for fractional Brownian motion and applications. Probab. Appl (N. Y.). Springer. --- 2008. --- 326 p. 

\bibitem{5}
Bogachev V. I. Gaussian measures. American Mathematical Society. Mathematical Surveys and Monographs. --- Volume 62. --- 1998. --- 433 p.

\bibitem{6}
Bondarenko V. V. Approximation of Time Series by Power Function of the Fractional Brownian Motion // J. Automation and Information Science. --- 2013. --- V. 45, I. 6. --- P. 82--86.

\bibitem{7}
Breton J. C., Nourdin I.  Error bounds on the non-normal approximation of Hermite power variations of fractional Brownian motion // Electronic Communications in Probability, 2008. --- V. 13. --- P. 482--493.

\bibitem{8}
Breton J. C., Nourdin I., Peccati G. Exact confidence intervals for the Hurst parameter of a fractional Brownian motion // Electronic Journal of Statistics, 2009. --- V. 3. --- P. 416--425.

\bibitem{9}
Coeurjolly J.-F. Simulation and identification of the fractional Brownian      motion: A bibliographical and comparative study / Coeurjolly J.-F. // Journal of statistical software. --- V. 5. --- Issue 7. --- 2000. --- P. 1--52.

\bibitem{10}
Coeurjolly J.-F. Estimating the parameters of a fractional Brownian motion by discrete variations of its sample paths / J.-F. Coeurjolly // Statistical Inference for Stochastic Processes, 2001. --- Vol. 4. --- No. 2. --- P. 199--227.

\bibitem{11}
Coeurjolly J.-F. Hurst exponent estimation of locally self-similar Gaussian processes using sample quantiles. The Annals of statistics / Coeurjolly J.-F. --- Vol. 36. --- 2008. --- 3. --- P. 1404--1434.

\bibitem{12}
Hu Y., Nualart D. Parameter estimation for fractional Ornstein-Uhlenbeck processes, Stat. Probab. Lett. 80 (2010). P. 1030--1038.

\bibitem{13}
Kubilius K., Mishura Y., Ralchenko K., Seleznjev O. Consistency of the drift parameter estimator for the discretized fractional Ornstein-Uhlenbeck process with Hurst index   $H \in \left(0, \frac{1}{2}\right)$. Electron. J. Stat. 9 (2015). P. 1799--1825.

\bibitem{14}
Lei P.,  Nualart  D.  A decomposition of the bifractional Brownian motion and some applications.  Statistics \& Probability Letters, 2009, V. 79, Issue 5, P. 619--624.

\bibitem{15}
Mandelbrot B. B. Une classe de processus stochastiques homothetiques a
soi: application a la loi climatologique de H. E. Hurst // Comptes Rendus de
l'Academie des Sciences. Paris. 1965. V. 240. P. 3274--3277.

\bibitem{16}
Mandelbrot B. B., van Ness J. W. Fractional Brownian motions, fractional noises and applications // SLAM Review. 1968. V. 10. 4. P. 422--437.

\bibitem{17}
Mandelbrot B. B. Fractals: Form, Chance, and Dimension. San Francisco: Freeman, 1977.

\bibitem{18}
Mandelbrot B. B. The Fractal Geometry of Nature / B. Mandelbrot. --- Freeman and Co., San Francisco, 1982. --- V. 89, issue 2. --- P. 460.

\bibitem{19}
Mishura Y. Stochastic Calculus for Fractional Brownian Motion and Related Processes. Lecture Notes in Mathematics. V. 1929 / Mishura Y. // Springer-Verlag, 2008. --- 392 p.

\bibitem{20}
Nourdin I. Asymptotic behavior of weighted quadratic and cubic variations of fractional Brownian motion / Nourdin I. // Ann. Probab. 36. --- 2008. --- Number 6. --- P. 2159--2175.

\bibitem{21}
Nourdin I.,  R\'{e}veillac A.    Asymptotic behavior of weighted quadratic variations of fractional Brownian motion: the critical case $H= \frac{1}{4}$ //  The Annals of Probability. --- 2009. --- V. 37, issue 6, P. 2200--2230.

\bibitem{22}
Nourdin I. Noncentral convergence of multiple integrals / I. Nourdin // Ann. Probab. --- Volume 37. --- 2009. --- Number 4. --- P. 1412--1426.

\bibitem{23}
Nourdin I. Density formula and concentration inequalities with Malliavin calculus / I. Nourdin, F. G. Viens // Electron. J. Probab. --- 14. --- 2009. --- P. 2287--2300.

\bibitem{24}
Nourdin I. Central and non-central limit theorems for weighted power variations of fractional Brownian motion / I. Nourdin, D. Nualart, C. Tudor. Ann. Inst H. Poincar\'{e} Probab Statist. --- V. 46. --- 2010. --- 4. --- P. 1055--1079.

\bibitem{25}
Nourdin I. Selected Aspects of Fractional Brownian Motion. Bocconi and  Springer Series-Verlag Italia. --- 2012. --- 122 p.     

\bibitem{26}
Nourdin I., Zintout R. Cross-variation of young integral with respect to long-memory fractional brownian motions. Probability   and    mathematical  statistics  Vol. 36, Fasc. 35 (2016), pp. 35--46.

\bibitem{27}
Nualart  D.    Fractional Brownian motion: stochastic calculus and applications.
International Congress of Mathematicians, 2006,  P. 1541--1562.

\bibitem{28}
Nualart D.,  Saussereau B.  Malliavin calculus for stochastic differential equations driven by a fractional Brownian motion. 
--- Stochastic Processes and their Applications (2009) V. 119, Issue 2, P. 391--409. 

\bibitem{29}
Peltier R. F. A new method for estimating the parameter of fractional Brownian motion / Peltier R. F., Levy Vehel J. // Rapport de recherch\'{e} de l'INRIA, 1994. --- 27 p. --- 2396.
\end{thebibliography}
\end{document}